\input amstex
\documentstyle{amams} 
\document
\annalsline{157}{2003}
\received{November 27, 2001}
\startingpage{939}
\font\tenrm=cmr10
\input amssym.def
\input amssym.tex

\catcode`\@=11
\font\twelvemsb=msbm10 scaled 1100

\font\ninemsb=msbm10 scaled 800
\newfam\msbfam
\textfont\msbfam=\twelvemsb  \scriptfont\msbfam=\ninemsb
  \scriptscriptfont\msbfam=\ninemsb
\def\msb@{\hexnumber@\msbfam}
\def\Bbb{\relax\ifmmode\let\next\Bbb@\else
 \def\next{\errmessage{Use \string\Bbb\space only in math
mode}}\fi\next}
\def\Bbb@#1{{\Bbb@@{#1}}}
\def\Bbb@@#1{\fam\msbfam#1}
\catcode`\@=12

 \catcode`\@=11
\font\twelveeuf=eufm10 scaled 1100
\font\teneuf=eufm10
\font\nineeuf=eufm7 scaled 1100
\newfam\euffam
\textfont\euffam=\twelveeuf  \scriptfont\euffam=\teneuf
  \scriptscriptfont\euffam=\nineeuf
\def\euf@{\hexnumber@\euffam}
\def\frak{\relax\ifmmode\let\next\frak@\else
 \def\next{\errmessage{Use \string\frak\space only in math
mode}}\fi\next}
\def\frak@#1{{\frak@@{#1}}}
\def\frak@@#1{\fam\euffam#1}
\catcode`\@=12

\def\varep{\varepsilon}
\def\eln{{\rm ln}\, }
\title{The Erd\H os-Szemer\'edi problem\\ on sum set and product set}  
\shorttitle{The Erd\H os-Szemer\'edi problem} 
\acknowledgements{Partially supported by NSA.}
 \author{Mei-Chu Chang}
   \institutions{University of California,
Riverside, CA \\ 
{\eightpoint {\it E-mail address\/}: mcc\@math.ucr.edu}}

 \vglue2pt \centerline{\bf  Summary}
\vglue12pt

The basic theme of this paper is the fact that if $A$ is a finite set of 
integers, then the sum and product sets cannot both be small.
  A precise formulation of this fact is Conjecture 1 below due to 
Erd\H os-Szemer\'edi [E-S]. (see also [El], [T], and
[K-T] for related aspects.) Only much weaker results or very special
 cases of this conjecture are presently known. One approach consists of 
assuming the sum set $A + A$ small and then deriving that the product set
 $AA$ is large (using Freiman's structure theorem) (cf.\ [N-T], [Na3]).  
We follow the reverse route and prove that if $|AA| < c|A|$,
 then $|A+A| > c^\prime |A|^2$ (see Theorem 1).  A quantitative version of
 this phenomenon combined with the Pl\"unnecke type of inequality (due to
Ruzsa) permit us to settle completely a related conjecture in [E-S] on the growth in $k$.  If
$$
g(k) \equiv \text{\rm min}\{|A[1]| + |A\{1\}|\}
$$
over all sets $A\subset \Bbb Z$ of cardinality $|A| = k$ and where $A[1]$
 (respectively, $A\{1\}$) refers to the simple sum (resp., product) of elements of $A$. (See (0.6), (0.7).)  It was conjectured in [E-S] that $g(k)$ grows faster than any power of $k$ for $k\rightarrow\infty$. 
We will prove here that $\eln g(k)\sim\frac{(\eln k)^2}{\eln \eln k}$ (see Theorem 2) which is the main
result of this paper.

\vglue16pt \centerline{\bf Introduction}
\vglue16pt

Let $A,B$ be finite sets of an abelian group.
The {\it sum set\/} of $A,B$ is
$$
A + B \equiv \{a + b \mid a \in A, b\in B\}. \tag 0.1
$$
We denote by
$$
hA\equiv A+\cdots+A\;\;\text{($h$ fold)}\tag 0.2
$$
the $h$-fold sum of $A$.

Similarly we can define the {\it product set} of $A,B$ and $h$-fold product of $A$.
$$
\align
AB &\equiv \{ab \mid a \in A, b \in B\}, \tag 0.3\\
A^h &\equiv A\cdots A\;\; \text{($h$ fold)}. \tag 0.4
\endalign
$$
If $B = \{b\}$, a singleton, we denote $AB$ by $b\cdot A$.
\medskip

In 1983, Erd\H os and Szemer\'edi [E-S] conjectured that for subsets of integers, the sum set and the product set cannot both be small.  Precisely, they made the following conjecture.
 
\nonumproclaim{Conjecture 1 {\rm (Erd\H os-Szemer\'edi)}}  For any $\varep > 0$ and any 
$h\in\Bbb N$ there is $k_0 =
k_0(\varep)$ such that for any $A \subset \Bbb N$ with $|A|\geq k_0${\rm ,}  
$$
|hA\cup A^h| \gg |A|^{h-\varep}. \tag 0.5
$$
\endproclaim

We note that there is an obvious upper bound $|hA\cup A^h|  \leq  2\! \left(\!\matrix |A| + h - 1\\h
\endmatrix\right)\!$.
\medskip
Another related conjecture requires the following notation of {\it simple sum} and {\it simple product}.
$$
\align
A[1] &\equiv \left\{\operatornamewithlimits{\dsize\sum}^k_{i=1} \varep_i a_i \mid a_i \in A, \varep_i = 0\;\;
\text{or}\;\; 1\right\},\tag 0.6\\
A\{1\} &\equiv \left\{\operatornamewithlimits{\dsize\prod}^k_{i=1} a^{\varep_i}_i \mid a_i \in A, \varep_i = 0\;\;\text{or}\;\;1\right\}. \tag 0.7
\endalign
$$

For the rest of the introduction, we only consider $A\subset \Bbb N$.
\medskip

\nonumproclaim{Conjecture 2   {\rm (Erd\H os-Szemer\'edi)}} 
 Let $g(k) \equiv \text{\rm min}_{|A| = k} \{|A[1]| + |A\{1\}|\}$. 
Then for any $t${\rm ,} there is $k_0 = k_0(t)$ such that for any $k\geq k_0, g(k) > k^t$.
\endproclaim

Toward Conjecture 1, all work has been done so far, are   for the case $h = 2$.
\medskip
Erd\H os and Szemer\'edi [E-S] got the first bound:
\medskip
\nonumproclaim{Theorem {\rm (Erd\H os-Szemer\'edi)}}  Let $f(k)\equiv\text{\rm min}_{|A|=k} |2A\cup A^2|$. 
 Then there are constants $c_1,c_2${\rm ,} such that
$$
k^{1+c_1} < f(k) < k^2 e^{-c_2 \frac{\eln  k}{\eln  \eln  k}}. \tag 0.8
$$
\endproclaim

Nathanson showed that $f(k) > ck^{\sssize \frac{32}{31}}$, with $c = 0.00028\ldots \  $.
 
At this point, the best bound is
$$
|2A\cup A^2| > c|A|^{\sssize 5/4} \tag 0.9
$$
obtained by Elekes [El] using the Szemer\'edi-Trotter theorem on line-incidences in the plane (see [S-T]).
 
On the other hand, Nathanson and Tenenbaum [N-T] concluded something stronger by assuming the sum set is small.  They showed

\nonumproclaim{Theorem {\rm (Nathanson-Tenenbaum)}}  If
$$
|2A| \leq 3|A| - 4, \tag 0.10
$$
then
$$
|A^2| \geqq \left( \frac{|A|}{\eln |A|}\right)^2. \tag 0.11
$$
\endproclaim 

Very recently, Elekes and Ruzsa [El-R] again using the
Szemer\'edi-Trotter theorem, established the following general inequality.

\nonumproclaim{Theorem {\rm (Elekes-Ruzsa)}} If $A \subset \Bbb R$ is a finite set{\rm ,} then
$$
|A+A|^4\, |AA|\,\eln |A|>|A|^6.\tag 0.12
$$
\endproclaim
 
In particular, their result implies that if 
$$
|2A|< c|A|,\tag 0.13
$$
then
$$
|A^2| \geqq \frac{|A|^2}{c^\prime \eln |A|}.\tag 0.14
$$

For further result in this direction, see [C2].
 
Related to Conjecture 2, Erd\H os and Szemer\'edi [E-S] have an upper bound:
 
\nonumproclaim{Theorem {\rm (Erd\H os-Szemer\'edi)}} \hskip-4pt 
 Let $g(k) \equiv \text{\rm min}_{|A| = k} \{|A[1]| + |A\{1\}|\}.\!$   There is a constant $c$ such that
$$
g(k) < e^{c\frac{(\eln k)^2}{\eln \eln k}}.\tag 0.15
$$
\endproclaim

Our first theorem is to show that the $h$-fold sum is big, if the product is small.

\nonumproclaim{Theorem 1}  Let $A \subset\Bbb N$ be a finite set.  If $|A^2| < \alpha|A|${\rm ,} then
$$
|2A| > 36^{-\alpha} |A|^2, \tag 0.16
$$
and
$$
|hA| > c_h (\alpha) |A|^h. \tag 0.17
$$
Here
$$
c_h (\alpha) = (2h^2 - h)^{-h\alpha}. \tag 0.18
$$
\endproclaim

Our approach is to show that there is a constant $c$ such that
$$
\int \big| {\sum}_{m\in A} e^{2\pi imx} \big|^{2h} dx < c|A|^h \tag 0.19
$$
by applying an easy result of Freiman's theorem (see the paragraph after
Proposition 10)
to obtain
$$
A \subset P \equiv \left\{\frac{a}{b}
 \left(\frac{a_1}{b_1}\right)^{j_1}\cdots \left(\frac{a_s}{b_s}\right)^{j_s} \bigm| 0\leq j_i < \ell_i\right\} \tag 0.20
$$
and carefully analyzing the corresponding trigonometric polynomials (see\break Proposition 8). These are estimates
 in the spirit of Rudin [R]. The constant $c$ here depends, of course, on $s$ and $h$.

In order to have a good universal bound $c$, we introduce the concept of multiplicative dimension of a finite set of integers, and derive some basic properties of it (see Propositions 10 and 11). 
 We expect more applications coming out of it.

Another application of our method together with a Pl\"unnecke type of
 inequality (due to Ruzsa) gives a complete answer to Conjecture 2.

\nonumproclaim{Theorem 2} 
 Let $g(k) \equiv \text{\rm min}_{|A| = k} \left\{ |A[1]| + |A\{1\}|\right\}$.  Then there is $\varep > 0$ such that
$$
k^{(1+\varep)\frac{\eln k}{\eln \eln k}} > g(k) > k^{(\frac{1}{8}-\varep)\frac{\eln k}{\eln \eln k}}.\tag 0.21
$$
\endproclaim

{\it Remark} 2.1 (Ruzsa). The lower bound can be improved to 
$k^{(\frac{1}{2}-\varep)\frac{\eln k}{\eln \eln k}}$.
We will give more detail after the proof of Theorem 2.
\vglue12pt

Using a result of Laczkovich and Rusza, we obtain the following result
related to a conjecture in [E-S] on undirected graphs.

\nonumproclaim{Theorem 3} Let 
 $G\subset A\times A$ satisfy $|G|>\delta |A|^2$.
Denote the restricted sum and product sets by
$$
\align
A\overset G\to +A &= \{ a+a'| (a, a') \in G\}\tag 0.22\\
A\overset G\to \times A&= \{aa'| (a, a') \in G\}.\tag 0.23
\endalign
$$
If
$$
|A \overset G\to\times A| < c|A|,\tag 0.24
$$
then
$$
|A\overset G\to +A|  > C (\delta, c) |A|^2.\tag 0.25
$$
\endproclaim

The paper is organized as follows:
 
In Section 1, we prove Theorem 1 and introduce the concept of multiplicative dimension.
 In Section 2, we show the lower bound of Theorem 2 and Theorem 3.
 In Section 3, we repeat Erd\H os-Szemer\'edi's upper bound of Theorem~2.

\demo{Notation}  We denote by $\lfloor a\rfloor$ the greatest integer $\leq a$, and by $|A|$ the cardinality of a set
$A$.
\enddemo 

{\it Acknowledgement}.   The author would like to thank J.\ Bourgain for various
 advice, and I.\ Ruzsa and the referee for many helpful comments.
 \vglue19pt

\section{Proof of Theorem 1}

Let $A\subset\Bbb N$ be a finite set of positive integers, and let $\Gamma_{h,A}(n)$ be the number of representatives of $n$ by the sum of $h$ (ordered) elements in $A$, i.e.,
$$
\Gamma_{h,A} (n) \equiv \left\vert \{(a_1,\ldots,a_h) \mid \sum a_i = n, a_i \in A\}\right\vert. \tag 1.1
$$
The two standard lemmas below provide our starting point.
 
\nonumproclaim{Lemma 3}  Let $A\subset\Bbb N$ be finite and let $h\in \Bbb N$.  If there is a constant $c$ such that
$$
\sum_{n\in h A} \Gamma^2_{h,A} (n) < c |A|^h, \tag 1.2
$$
then
$$
|hA| > \frac{1}{c} |A|^h.  \tag 1.3
$$
\endproclaim

\demo{Proof}  Cauchy-Schwartz inequality and the hypothesis give
$$
\align
|A|^h = \sum_{n\in hA} \Gamma_{h,A} (n) &\leq |hA|^{\sssize 1/2} \left(\sum_{n\in hA} \Gamma^2_{h,A}
(n)\right)^{\sssize 1/2}\\
      &< |hA|^{\sssize 1/2} c^{\sssize{1/2}} |A|^{h/2}.\\
\noalign{\vskip-24pt}
\endalign
$$
\enddemo
\vglue8pt

\nonumproclaim{Lemma 4} The following equality holds\/{\rm :}\/
$$
\sum_{n\in hA} \Gamma^2_{h,A} (n) = \left(\|\sum_{m\in A} e^{2\pi imx} \|_{{}_{\sssize 2h}}\right)^{\sssize 2h}.
$$
\endproclaim

\demo{Proof}
$$
\align
\left(\left\|\sum_{m\in A} e^{2\pi imx}\right\|_{{}_{2h}}\right)^{\sssize 2h} &=
 \int \left\vert \sum_{m\in A}
e^{2\pi imx}\right\vert^{\sssize 2h}\; dx\\ \noalign{\vskip4pt} &=\
 \int \left\vert \left(\sum_{m\in A} e^{2\pi imx}\right)^h\right\vert^2\; dx\\ \noalign{\vskip4pt}
&= \int \left\vert\left(\sum_{n\in hA} \Gamma_{h,A} (n) e^{2\pi inx}\right)\right\vert^2\; dx\\ \noalign{\vskip4pt}
&= \sum_{n\in hA} \Gamma^2_{h,A} (n).
\endalign
$$
The last equality is Parseval equality.  
\enddemo

From Lemmas 3 and 4, it is clear that to prove Theorem 1, we want to find a constant $c$ such that
$$
\left(\left\|\sum_{m\in A} e^{2\pi imx}\right\|_{\sssize 2h}\right)^{\ssize{2}} < c |A|. \tag 1.4
$$
In fact, we will prove something more general   to be used in the inductive
argument.
 
\nonumproclaim{Proposition 5} Let $A \subset\Bbb N$ be a finite set with
 $|A^2| < \alpha
|A|$. Then for any $\{d_a \}_{a \in A} \subset\Bbb R_+${\rm ,}  
 $$
\left(\left\|\sum_{a\in A} d_a e^{2\pi iax}\right\|_{{}_{2h}}\right)^{\sssize 2} < c\, \sum d^2_a
\tag 1.5
$$
for some constant $c$ depending on $h$ and $\alpha$ only.
\endproclaim

For a precise constant $c$, see Proposition 9.

The following proposition takes care of the special case of (1.5) when there
 exists a prime $p$ such that for every nonnegative integer $j, p^j$ appears
in the prime factorization of at most one element in $A$.  It is also the
 initial step of our iteration.

First, for convenience, we use the following:
\demo{Notation}  We denote by $\langle G\rangle^+$, the set of linear combinations of elements in $G$ with
\pagebreak coefficients in $\Bbb R^+$.\enddemo

\nonumproclaim{Proposition 6}  Let $p$ be a fixed prime{\rm ,} and let
$$
F_j(x) \in \left\langle\left\{e^{2\pi ip^jnx}\bigm|  n\in\Bbb {N}, (n,p) = 1\right\}\right\rangle^+. \tag 1.6
$$
Then
$$
\left(\left\|\sum_j F_j\right\|_{{}_{2h}}\right)^{\sssize 2} \leq c_h \sum_j \left\|F_j\right\|^{\sssize
2}_{{}_{2h}},\;\;\text{where
$c_h = 2h^2 - h$}.
\tag 1.7
$$
\endproclaim

\demo{Proof}  To bound $\int |\sum_j F_j |^{2h} dx$, we expand $|\sum_jF_j|^{2h}$ as
$$
\left(\sum F_j\right)^h \left(\sum \overline{F}_j\right)^h.
\tag 1.8
$$
Let
$$
F_{j_1}\cdots F_{j_h}\;\overline{F}_{j_{h+1}}\cdots\overline{F}_{j_{2h}}
\tag 1.9
$$
be a term in the expansion of (1.8). After rearrangement, we may assume 
$j_1\leq\cdots\leq j_h$, and $j_{h+1} \leq\cdots\leq j_{2h}$.

When (1.9) is expressed as a linear combination of trignometric functions, a typical term is of the form
$$
ne^{2\pi ix (p^{j_1}n_1+\cdots+p^{j_h} n_h - p^{j_{h+1}}
n_{h+1}-\cdot\cdot-p^{j_{2h}}n_{2h})}.
\tag 1.10
$$

We note that the integral of (1.10) is $0$, if the expression in the parenthesis in (1.10) is nonzero.  In particular, independent of the $n_i$'s, the integral of (1.10) is $0$, if
$$
j_1\not= j_2\leq j_{h+1}, \quad \text{or}\quad j_1\not= j_{h+1}\leq\text{\rm min}\{j_2,j_{h+2}\},
\quad \text{or}\quad j_{h+1}
\not= j_{h+2} \leq j_1.
\tag 1.11
$$
Therefore, if any of the statements in (1.11) is true, then the integral of (1.9) is $0$.

\medskip

We now consider the integral of (1.9) where the index 
set $\{j_1,\ldots,j_{2h}\}$ does not satisfy any of the conditions in (1.11).  
For the case $j_1 = j_2 \leq j_{h+1}$, we see that in an ordered set of $h$ 
elements coming from the expansion of (1.8) (before the rearrangement), 
there are exactly $\binom h2$ choices for the positions 
of $j_1, j_2$.  On the other hand, if $F_{j_1}F_{j_2}$ is factored out, the 
rest is symmetric with respect to $j_3,\ldots,j_h$, and $j_{h+1},\ldots,j_{2h}$, i.e., all the terms involving $j \equiv j_1 = j_2 \leq j_{h+1}$ are simplified to
$$
\binom h2 \left(F_j\right)^2 \left(\sum_{k\geq j} F_k\right)^{h-2}.
\tag 1.12
$$
With the same reasoning for the other two cases, we conclude that

$$
\align
\noalign{\vskip-18pt}
\left(\Big\|\sum_j F_j\Big\|_{\sssize 2h}\right)^{2h} =& \binom h2 \sum_j \int F^2_j \left(\sum_{k\geq j}
F_k\right)^{h-2}
\left(\sum_{k\geq j}
\overline{F}_k\right)^h\, dx\\
 &+ h^2 \sum_j \int |F_j|^2 \left(\sum_{k\geq j} F_k \sum_{k\geq j}\overline{F}_k\right)^{h-1}\, dx\\
&+ \binom h2 \sum_j \int \overline{F}_j^2 \left(\sum_{k\geq j} F_k\right)^h \left(\sum_{k\geq j}
\overline{F}_k\right)^{h-2}\, dx.
\endalign
$$
The right-hand
side is
$$
\align
&\leq \left[h^2 + 2 \binom h2\right] \sum_j \int |F_j|^2 \left\vert\sum_{k\geq j} F_k\right\vert^{\sssize 2h-2}\,
dx\\ &\leq (2h^2 - h) \sum_j \|F^2_j\|_{{}_{h}} \|\left(\sum_{k\geq j} F_k\right)^{2h-2}\|_{\frac{h}{h-1}}\\
&= (2h^2 - h) \sum_j \|F_j\|^{\sssize 2}_{{}_{2h}} \left(\Big\|\sum_{k\geq j} F_k\Big\|_{{}_{2h}}\right)^{\sssize
2h-2}. 
\endalign
$$
The last inequality is H\"older inequality.
 
Now, the next lemma concludes the proof of Proposition 6. 
\enddemo
 
\nonumproclaim{Lemma 7}  Let $F_k\in\langle\{e^{2\pi i m_k x} \mid  m_k\in \Bbb {Z}\}\rangle^+$.  Then
$$
\left\|\sum_k F_k\right\|_{{}_{2h}}\geq \left\|\sum_{k\geq j} F_k\right\|_{{}_{2h}},\quad\text{for any}\quad j.
\tag 1.13
$$
\endproclaim

{\it Proof}.
$$
\align
&\int \left\vert\sum_k F_k\right\vert^{2h} dx \\
&\qquad =  \int \left(\sum_{k\geq j} F_k + \sum_{k < j} F_k\right)\\
&\qquad\quad \cdots \left(\sum_{k\geq j} F_k + \sum_{k < j}
F_k\right)\left(\sum_{k\geq j} \overline{F}_k + \sum_{k<j}\overline{F}_k\right)\cdots\left(\sum_{k\geq j} 
\overline{F}_k
+
\sum_{k < j}
\overline{F}_k\right)\, dx \\
&\qquad  \geq \int \left(\sum_{k\geq j} F_k \sum_{k\geq j} \overline{F}_k\right)^h dx =\left(\|\sum_{k\geq j}
F_k\|_{\sssize 2h}\right)^{2h}. \\
\endalign
$$
The inequality holds because the coefficients of the trignometric functions (as in (1.10)) in the expansion are all
positive. 
\hfill\qed\vglue12pt
 
{\it Remark} 7.1.  This is a special case of a general theorem in martingale theory.
 
\nonumproclaim{Proposition 8}  Let $p_1,\cdots, p_t$ be distinct primes{\rm ,} and let
$$
F_{j_1,\ldots,j_t} (x) \in \left\langle \left\{e^{2\pi i p_1^{j_1}\cdots p^{j_t}_t nx} \bigm| n\in \Bbb {N}, (n, p_1\cdots
p_t) = 1\right\}\right\rangle^+.
\tag 1.14
$$
Then
$$
\left\|\sum_{j_1,\ldots,j_t} F_{j_1,\ldots,j_t}\right\|^{\sssize 2}_{\sssize{2h}} \leq c^t_h \sum_{j_1,\ldots,j_t}
\|F_{j_1,\ldots,j_t}\|^{\sssize 2}_{\sssize{2h}}, \;\;\text{where $c_h = 2h^2 - h$}. \tag 1.15
$$
\endproclaim  

\demo{Proof} 
We do induction on $t$. 
The left-hand side of (1.15) becomes
$$
\left\|\sum_{j_1}\sum_{j_2,\ldots,j_t} F_{j_1,\ldots,j_t}\right\|^{\sssize 2} \leq
c_h \sum_{j_1} \left\|\sum_{j_2,\ldots,j_t} F_{j_1,\ldots,j_t}\right\|^{\sssize 2} \leq
c_h \sum_{j_1}c^{t-1}_h\sum_{j_2,\ldots,j_t}\| F_{j_1,\ldots,j_t}\|^{\sssize 2},
$$ 
which is the right-hand side.  
\enddemo

Proposition 5 is proved, if we can find a {\it small} $t$ such that the Fourier
transform of $F_{j_1,\ldots,j_t}$ is supported at one point and such $t$ is
bounded by $\alpha$. So we introduce the following notion.

\demo{Definition}  Let $A$ be a finite set of positive rational
 numbers in lowest terms (cf.\ (0.20)). Let $q_1,\ldots,q_\ell$ be all the prime
factors in the obvious prime factorization of elements in $A$.  For $a \in A$,
 let $a = q^{j_1}_j\cdots q^{j_\ell}_\ell$ be the prime factorization of $a$.
Then the map $\nu : A\rightarrow \Bbb {R}^\ell$ by sending $a$ to $(j_1,\ldots,
j_\ell)$ is one-to-one.  The {\it multiplicative dimension} of $A$ is the
 dimension of the smallest (affine) linear space in $\Bbb {R}^\ell$
 containing $\nu(A)$.
\enddemo

We note that for any nonzero rational number $q$, $q \cdot A$ and $A$ have the same multiplicative dimension, since $\nu (q\cdot A)$ is a translation of $\nu(A)$.
 
The following proposition is a more precise version of Lemma 5.
\vfill
\nonumproclaim{Proposition 9}  Let $A\subset\Bbb{N}$ be finite with {\rm mult.}$\dim (A) = m$.  Then
$$
\left(\left\|\sum_{a\in A} d_a e^{2\pi iax}\right\|_{\sssize 2h}\right)^{\sssize 2} < c_h^{\sssize{m}} \sum d^{\sssize
2}_a,\quad\text{ where }\quad c_h = 2h^{\sssize 2} - h.
\tag 1.16
$$
\endproclaim

\demo{Proof}  To use (1.15) in Proposition 8, we want to show that there are
primes $q_1,\ldots,q_m$ such that a term of the trigonometric
 polynomial in the 
left-hand side of (1.15), when expressed in terms of the notation in (1.14),
is $F_{j_1,\ldots,j_m}= d_a e^{2\pi i q_1^{j_1}\cdots q^{j_m}_m nx}$.
In other words, we want to show that among the prime factors $q_1,\ldots,q_\ell$ of elements in $A$, there are $m$ of them, say $q_1,\ldots,q_m$ such that

$(\ast)$
 for all $(j_1,\ldots,j_m) \in\Bbb {Z}^m$, there is at most one $a\in A$ such that $q_1^{j_1}\cdots
q^{j_m}_m$  is part of the prime factorization of $a$.
This is equivalent to
 
$(\ast\ast)$
$\pi \circ \nu$ is injective, where $\nu$ is as in the definition of multiplicative dimension and $\pi :
\Bbb {R}^\ell \rightarrow \Bbb {R}^m$ is the projection to the first $m$ coordinates.

Since $\dim\nu(A) = m$, $(\ast\ast)$ is clear after some permutation of the $q_i$'s.  
\enddemo

\nonumproclaim{Proposition 10}  Let $A\subset\Bbb {N}$ be finite with
 {\rm mult.}$\dim A = m$.
Then 
$$
\sum_{n\in hA}\Gamma^2_{h,A} (n) < c^{{}^{mh}}_h |A|^h,\;\;\text{ where  }\;\; c_h = 2h^2 - h.
\tag 1.17
$$
\endproclaim

\demo{Proof}  This is a consequence of Lemma 4 and Proposition 9 (with\break $d_a = 1$). 
\enddemo

The hypothesis of Theorem 1 gives a universal
 bound on the multiplicative
 dimension of $A$ by applying Freiman's 
theorem (cf. [Fr1], [Fr2], [Fr3], [Bi], [C1], [Na1]). In fact, we do not 
need the full content of {\it the} Freiman's theorem, but a much easier result by Freiman.
  A small modification (over $\Bbb Q$ instead of over~$\Bbb R$) of Lemma 4.3
 in [Bi] is sufficient. (As Ruzsa pointed out   it is also Lemma 1.14
 in [Fr1].)
\medskip

\nonumproclaim{Theorem {\rm (Freiman)}}  Let $G\subset\Bbb R$ be a subgroup and
 $A_1 \subset G$ be finite.  If there is a constant
 $\alpha${\rm ,} $\alpha < \sqrt{|A_1|}${\rm ,} such that $|2A_1| < \alpha |A_1|${\rm ,}
 then there is an integer
$$
s \leq \alpha
$$
such that $A_1$ is contained in an $s$\/{\rm -}\/dimensional proper progression
 $P_1${\rm ;} i.e.{\rm ,} there exist $\beta,\alpha_1,\ldots,\alpha_s \in G$ and
 $J_1,\cdots,J_s \in \Bbb N$ such that
$$
A_1 \subset P_1 = \{\beta + j_1\alpha_1 +\cdots+ j_s \alpha_s \mid 0\leq j_i < J_i\},
$$
and $|P_1| = J_1\cdots J_s$.
\endproclaim

Note that if $|A_1| > \frac{\lfloor \alpha\rfloor\lfloor\alpha + 1\rfloor}{2(\lfloor\alpha+1\rfloor - \alpha)}$,
 then $s\leq \lfloor\alpha - 1\rfloor$.
\vglue4pt
 
Recall that the full Freiman  theorem also permits one to state a bound
$J_1 \cdots J_s < c(\alpha)|A_1|$. However this additional information will
not be used in what follows.

We would like to work on a sum set instead of a product set.  So we define
$$
A_1 \equiv \eln  A = \{\eln  a\mid a\in A\}.
\tag 1.18
$$
Note that ${\rm ln}$ is an isomorphism between the two groups $(\Bbb{Q}^+,\, \cdot\, )$ and $(\eln  \Bbb{Q}^+,   + )$.
 
Applying the theorem to $A_1\subset \eln  \Bbb {Q}^+$, then pushing back by $({\rm ln})^{-1}$, we have
$$
A \subset P\equiv \left\{\frac{a}{b} (\frac{a_1}{b_1})^{j_1}\cdots(\frac{a_s}{b_s})^{j_s}\mid 0 \leq j_i < J_i\right\}\subset \Bbb {Q}^+,
\tag 1.19
$$
where $a,b,a_i,b_i,J_i \in \Bbb {N}$, and $(a,b) = 1, (a_i,b_i) = 1$.  Moreover, $s\leq \lfloor\alpha - 1\rfloor$ and different ordered sets $(j_1,\cdots,j_s)$ represent different rational numbers.  Clearly,
$$
\text{mult.}\dim A \leq \text{mult.dim }P = \dim E \leq s\leq \lfloor\alpha - 1\rfloor,
\tag 1.20
$$
where $E$ is the vector space generated by $\nu (\frac{a_1}{b_1}),\cdots,\nu (\frac{a_s}{b_s})$.

Therefore, we have
 
\nonumproclaim{Proposition 11}  Let $A\subset\Bbb{N}$ be a finite set. 
If $|A|^2 < \alpha |A|$ for some constant~$\alpha${\rm ,}  $\alpha < |A|^{1/2}${\rm ,} 
then
 {\rm mult.}$\dim A\leq \alpha$.
  Furthermore{\rm ,} if $|A| > \frac{\lfloor\alpha\rfloor\lfloor\alpha+1\rfloor}
{2(\lfloor\alpha+1\rfloor-\alpha)}${\rm ,} then {\rm mult.}$\dim$ $A\leq
 \lfloor\alpha-1 \rfloor$.
\endproclaim

Putting Propositions 10 and 11 together, we have
  
\nonumproclaim{Proposition 12}\hskip-3pt  Let $A\subset\Bbb N$ be finite.  
If $|A^2| < \alpha |A|$ for some constant~$\alpha${\rm ,} $\alpha < |A|^{1/2}${\rm ,} then
$$
\sum_{n\in hA} \Gamma^2_{h,A} (n) < c^{\alpha h}_h |A|^h,\quad\text{where }\quad c_h = 2h^2 - h.
$$
\endproclaim

Now, Theorem 1 follows from Proposition 12 and Lemma 3.\hfill\qed

\section{Simple sums and products}
 
In this section we will prove the lower bound in Theorem 2.
 
Let $A\subset\Bbb N$ be finite.  We define
$$
g(A)\equiv |A[1]| + |A\{1\}|,
\tag 2.1
$$
where $A[1]$ and $A\{1\}$ are the simple sum and simple product of $A$.  (See (0.6), (0.7) for precise definitions.)
\medskip

We will show that for any $\varep$ and any $A\subset\Bbb N$ with $|A| = k\gg 0$,
$$
g(A) > k^{(\frac18 - \varep)\frac{\eln k}{\eln \eln k}}.
\tag 2.2
$$
For those who like precise bounds, we show:

For $0 < \varep_1,\varep_2 < \frac12$,
$$
g(A) > e^{-3} \lfloor k^{\frac12 - \varep_2} \rfloor^{\lfloor(\frac14 - \frac
{\varep_1}{2})\frac{\eln k}{\eln \eln k}\rfloor},
\tag 2.3
$$
if $|A| = k$ is large enough such that
$$
\eln \eln k > \frac{\sqrt{2}}{8\varep_1},
\tag 2.4
$$
and
$$
\frac{\eln k}{\eln \eln k} > \frac{2}{\varep_2}.
\tag 2.5
$$
 
\nonumproclaim{Proposition 13}  Let $B\subset\Bbb N$ be finite with 
{\rm mult.}$\dim B = m$.  Then for any $h_1 \in\Bbb N${\rm ,}  
$$
|h_1 B\cap B[1]| > \left[\frac{|B|}{(2h^2_1 - h_1)^{m+1}}\right]^{h_1}.
\tag 2.6
$$
\endproclaim
 
\demo{Proof}  Since $h_1 B\cap B[1]$ is the set of simple sums with exactly $h_1$ summands, we have
$$
\pmatrix |B|\\h_1\endpmatrix \leq \sum_{n\in h_1B\cap B[1]} \Gamma _{h_1,B} (n). \tag 2.7
$$
Therefore,
$$
\align
\left(\frac{|B|}{h_1}\right)^{h_1} &< (h_1 B\cap B[1])^{1/2} \left(\sum_{n\in h_1B} \Gamma_{h_1,B}^2
(n)\right)^{1/2}\\ &\leq (h_1 B\cap B[1])^{1/2}\; [(2 h^2_1 - h_1)^{mh_1}|B|^{h_1}]^{1/2}.
\endalign
$$
The first inequality is because of  the Cauchy-Schwartz inequality and the fact that $h_1B\cap B[1]\subset h_1B$. 
The second inequality is Proposition 10.  \enddemo

{\it Remark} 13.1.  Clearly, from our proof, the denominator in (2.6) can be replaced
by  $ (2h_1 ^2 - h_1)^m h_1 ^2$.
 
\nonumproclaim{Proposition 14}  Let $B\subset\Bbb N$ with $|B| \geq \sqrt{k}$ and
{\rm mult.}$\dim B = m$.  For any $0 < \varep_1<\frac12${\rm ,} if
$$
m+1 \leq \left(\frac14 - \frac{\varep_1}{2}\right) \frac{\eln k}{\eln \eln k}, \tag 2.8
$$
then
$$
g(B) > k^{\varep_1\lfloor\frac{\eln k}{\sqrt{2}}\rfloor}.\tag 2.9
$$
\endproclaim

\demo{Proof}  Inequality (2.8) is equivalent to
$$
(\eln k)^{2m+2} \leq k^{{\sssize{1/2}} - \varep_1}. \tag 2.10
$$
In Proposition 13,  we take  $h_1 = \lfloor\frac{\eln k}{\sqrt{2}}\rfloor$ .  This gives
$$
2 h_1^2\leq (\eln k)^2.\tag 2.11
$$
Combining (2.11), (2.10) and (2.6), we have
\vglue12pt
\hfill ${\displaystyle
g(B) > \left|h_1B\cap B[1]\right| > \left(\frac{k^{\sssize{1/2}}}{k^{\sssize{1/2} - \varep_1}}\right)^{\lfloor\frac{\eln
k}{\sqrt{2}}\rfloor} = k^{\varep_1\lfloor
\frac{\eln k}{\sqrt{2}}\rfloor}.
}$
\enddemo

{\it Remark} 14.1.
  Let $A\subset\Bbb N$ with $|A| = k, k\gg 0$ (see (2.4)).  The set $B$ in Proposition 14 will be taken as a subset of $A$.  Then the bound in (2.9) is bigger than that in (2.2), and our proof is done.  Therefore for the rest of the section, we assume
$$
\text{mult.dim}B\geq \lfloor(\frac14 - \frac{\varep_2}{2}) \frac{\eln k}{\eln \eln k}\rfloor, \text{ for any $B\subset A$ with $|B| > \sqrt{k}$}.
\tag 2.12
$$
\vglue12pt

We need the following:
\demo{Notation}   We denote $B^\prime\equiv\nu(B)$ for any $B\subset A$, where $\nu= A\rightarrow\Bbb Z^\ell$ is as
in the definition of multiplicative dimension.
\enddemo
 
Note that
$$
|B^\prime [1]| = |B\{1\}|.\tag 2.13
$$

We will use the following Pl\"unnecke type of inequality due to Ruzsa.

\nonumproclaim{Ruzsa's Inequality {\rm [Ru2]}}  For any $h, \ell\in\Bbb N${\rm :}
$$
\text{If } |M+N| \leq \rho|M|, \text{ then } |hN - \ell N|\leq \rho^{h+\ell} |M|.
$$
\endproclaim

\demo{Proof of {\rm (2.2)}}  We divide $A$ into $\lfloor\sqrt{k}\rfloor$ pieces $B_1,B_2,\cdots,$ each of cardinality at least $\sqrt{k}$.  For $0 < \varep_2 < \frac12$, let
$$
\rho = 1 + k^{-1/2 + \varep_2}, \tag 2.14
$$
and let
$$
A_s \equiv \bigcup^s_{i=1} B_i. \tag 2.15
$$
There are two cases:
\medskip
\noindent
(i)\quad For all $s, |(A_s\bigcup B_{s+1})^\prime [1]| > \rho|A^\prime_s [1]|$.
Iterating gives
$$
|A^\prime[1]| = |(B_1\cup B_2\cup\cdots)^\prime[1]| > \rho^{\sqrt{k}-2} \sqrt{k}.\tag 2.16
$$
Therefore
$$
\align
 g(A)  &> |A\{1\}| = |A^\prime[1]| > e^{(\sqrt{k} - 2)\eln \rho+\frac12 \eln k}\tag 2.17
\\
& > e^{(\sqrt{k}-2)\frac45 k^{-1/2+\varep_2} + \frac12 \eln k}\\
& > e^{\frac45 k^{{}^{\varep_2}}}.
\endalign
$$
Inequality (2.5) is equivalent to
$$
k^{\varep_2} > (\eln k)^2.
$$
which is certainly stronger than what we need to show (2.2).

\vglue4pt
(ii)\quad  There exists $s$ such that $|(A_s\cup B_{s+1})^\prime[1]|\leq \rho|A^\prime_s[1]|$.
We use the fact that $(A_s\cup B_{s+1})^\prime[1] = A^\prime_s[1] + B'_{s+1} 
[1]$, and Ruzsa's inequality (with $h = h_2 + 1,\ell = 1$) to obtain
$$
|(h_2 + 1)B^\prime_{s+1}[1] - B^\prime_{s+1}[1]|\leq \rho^{h_2 + 2} |A^\prime_s[1]|. \tag 2.18
$$
Let $m = \text{mult.dim}B_{s+1}$.  For a set $B$, for $h\in\Bbb N$, denote
$$
B[h]\equiv\left\{\sum \varep_i x_i\mid \varep_i = 0,\ldots,h, x_i\in B\right\}. \tag 2.19
$$
The left-hand side of (2.18) is
$$
\align
&\geq |h_2B^\prime_{s+1}[1]|\tag 2.20
\\
&\geq |B^\prime_{s+1} [h_2]|\\
&\geq h^m_2.
\endalign
$$

We take $h_2 = \lfloor k^{\sssize{1/2} - \varep_2}\rfloor$.  Then in the 
right-hand
side of (2.18),
$$
\align
\rho^{h_2+2}    &\leq (1 + k^{- 1/2 + \varep_2})^{k^{1/2-\varep_2} + 2}\tag 2.21
\\
&< (e^{k^{-1/2+\varep_2}})^{k^{1/2-\varep_2}+2}\\
&< e^3.
\endalign
$$
Therefore, (2.18), (2.20) and (2.21) imply
$$
\align
g(A) &> g(A_s)\\
&> |A^\prime_s[1]| \\
&> e^{-3} h^m_2\\
  &> e^{-3} \lfloor k^{1/2 - \varep_2}\rfloor^{\lfloor (\frac14 - \frac{\varep_1}{2})\frac{\eln k}{\eln \eln k}\rfloor}.
\endalign
$$
The last inequality follows from our choice of $h_2$ and Remark 14.1.  
\enddemo
 
\demo{Proof of Remark {\rm 2.1}} In Proposition 14, if we take $B$ with $|B| \geq \frac {k}{2}$,
then we will replace (2.8), and (2.9) by
$$
m+1 \leq \frac12 \left(1 - \varep_1\right) \frac{\eln k}{\eln \eln k}, \tag $2.8'$ 
$$
and
$$
g(B) > \left(\frac{k^{\varep_1}}{2}\right)^{\lfloor\frac{\eln k}{\sqrt{2}}\rfloor}.\tag $2.9'$
$$
Let
$$
m_0 = \lfloor \frac12 \left(1 - \varep_1\right) \frac{\eln k}{\eln \eln k} \rfloor.
$$
Then (2.12) will be replaced by
$$
\text{mult.dim}B\geq m_0, 
\text{ for any $B\subset A$ with $|B| > \frac{k}{2}$}.
\tag $2.12'$
$$
Now we modify the proof of (2.2).

 Since $|A|=k> \frac{k}{2}$,
we have mult.dimA $\geq m_0$. So there is $B_1 \subset A$ with
mult.dimB$_1 =m_0$ and $|B_1|= m_0 +1$. Similarly, we have
$B_2 \subset A - B_1$ with mult.dimB$_2 =m_0$ and $|B_2|=m_0 +1$.
We continue this process until $r > \frac{k}{2(m_0 +1)}$. We have
$$
A \supset B_1\cup \cdots \cup B_r
$$
with
$$
\text{mult.dim}B_i = m_0, 
$$
and 
$$
|B_i|= m_0 +1.
$$
With more replacements,
$$
\frac{\eln k}{\eln \eln k} > \frac{3}{\varep_2}.\tag $2.5'$
$$
and
$$
\rho = 1 + k^{-1 + \varep_2}, \tag $2.14'$
$$
Identical arguments give 
\vglue12pt 
\hfill ${\displaystyle
g(A) > e^{-3} \lfloor k^{1 - \varep_2}\rfloor^{\lfloor (\frac12 - \frac{\varep_1
}{2})\frac{\eln k}{\eln \eln k}\rfloor}. 
}$
\enddemo 

\vglue8pt
\demo{Sketch of Proof of Theorem $3$} Let $|A|=N$. Then the  Laczkovich-Ruzsa
 theorem [L-R] and
(0.24) give $A_1 \subset A$ with
$$
|A_1A_1| < c^{\prime} N, \tag 2.22
$$
and
$$
|G \cap (A_1 \times A_1)| > \delta' N^2. \tag 2.23
$$
The weak Freiman theorem and (2.22) imply
$$
\text{mult.dim}A _1 < c^{\prime}.\tag 2.24
$$
It follows from Proposition 10 (with $h=2$) and   the proof of Lemma 4  that
$$
\beta\equiv |\{ (n_1,n_2,n_3,n_4) \in A_1 ^4 | n_1 -n_2 +n_3 -n_4 =0 \}|<36 ^{c'}N^2.\tag 2.25
$$
Hence
$$
\delta ' N^2 < \sum_{n \in A_1 \overset G\to +A_1} |\{(n_1,n_2) \in A_1 ^2 | n=n_1 +n_2\}| <
|A_1\overset G\to +A_1|^{\frac{1}{2}} \beta ^{\frac{1}{2}} .\tag 2.26
$$
The first inequality is (2.23), while the second one is 
the Cauchy-Schwartz inequality.

Therefore, (2.25) and (2.26) give
\vglue12pt \hfill  ${\displaystyle
|A\overset G\to +A| \geq |A_1\overset G\to +A_1| \geq \frac{(\delta ')^2 N^4}{\beta} >C N^2.  
}$\hfill\qed
\enddemo

\section{The example}

In this section for completeness we repeat a family of examples by Erd\H os-Szemer\'edi which provide the upper bound in Theorem 2.  Precisely, we will show

\nonumproclaim{Proposition 15}  Given $\varep_3 > 0${\rm ,} for $J$ so large that
$$
\frac{\eln J}{\eln \eln J} > \frac{1}{\varep_3}, \tag 3.1
$$
there is a set $A$ of cardinality $|A| = k\equiv J^J${\rm ,}  such that
$$
g(A) < 2 k^{(1 + \varep)\frac{\eln k}{\eln \eln k}}, \tag 3.2
$$
where
$$
\varep = 3\varep_3 + \varep^2_3. \tag 3.3
$$
\endproclaim

The example really comes from the proof of the lower bound of Theorem~2.
 
Let $p_1,\cdots, p_J$ be the first $J$ primes, and let
$$
A\equiv \left\{p_1^{{}^{j_1}}\cdots p_J^{{}^{j_J}} \mid 0 \leq j_i < J\right\}. \tag 3.4
$$
Then
$$
k\equiv |A| = J^J. \tag 3.5
$$

We will use the following relations between $k$ and $J$.
 
\nonumproclaim{Lemma 16}  Let $k,J$ be as in {\rm (3.5).}  Then  
\vglue4pt

{\rm (i)}  $\eln k = J \eln J$.
\vglue4pt
{\rm (ii)} $\eln \eln k = \eln J + \eln \eln J$.
\vglue4pt
\noindent
If $J$ and $\varep_3$ satisfy {\rm (3.1),} then
\vglue4pt
{\rm (iii)} $ \eln \eln k < (1 + \varep_3) \eln J$.
\vglue4pt
{\rm (iv)}  $J < (1 + \varep_3) \frac{\eln k}{\eln \eln k}$.
\vglue4pt
{\rm (v)}  $J^2 < (1 + \varep^\prime) (\frac{\eln k}{\eln \eln k})^2${\rm ,}
 where $\varep^\prime = 2\varep_3 + \varep^2_3$.
\endproclaim

\demo{Proof}  Each one follows immediately from the preceding one. 
For (iii) implying (iv), we use $J = \frac{\eln k}{\eln J}$.  
\enddemo

{\it Remark} 16.1.  The inequality $\frac{\eln J}{\eln \eln J} > \frac{1}{\varep_3}$ clearly implies
$$
\frac{\eln k}{\eln \eln k} > \frac{1}{\varep_3}.
\tag 3.6
$$

\nonumproclaim{Lemma 17} {\rm (i)} For all $a \in A${\rm ,} $a < (\eln k)^{J^2},$ \vglue4pt
{\rm (ii)}   $|A[1]| < k(\eln k)^{J^2},$
\vglue4pt
{\rm (iii)}  $|A\{1\}| < (kJ)^J$.
\endproclaim

\demo{Proof} 
(i)  For $a \in A$, (3.4) gives
$$
\align
a < \left(\operatornamewithlimits{\prod}_{i<J} p_i\right)^J &<
 \left(\operatornamewithlimits{\prod}_{i < J} i \eln i\right)^J\\
 &< \left(J^J (\eln J)^J\right)^J\\
&= \left(J \eln J\right)^{J^2}\\
&= (\eln k)^{J^2}. 
\endalign
$$
The second inequality is by the Prime Number Theorem.  
 The last equality is Lemma 16 (i).
\vglue4pt
(ii) follows from (i).
\vglue4pt
(iii)  We see that
$$
A\{1\} = \left\{p_1^{\sum^k_{s=1} j_1^{(s)}}\cdots p_J^{\sum^k_{s=1} j^{(s)}_J}  \mid 0 \leq j_i^{(s)} < J\right\}. \tag 3.7
$$
Since $\sum^k_{s=1} j^{(s)}_i < k J$, (iii) holds. \enddemo

\demo{Proof of Proposition {\rm 15}}  Lemma 17 (ii) and Lemma 16 (v) give 
$$
\align
|A[1] | &< k (\eln k)^{(1+\varep^\prime)(\frac{\eln k}{\eln \eln k})^2}\tag 3.8
\\
&= e^{\eln k + (1 + \varep^\prime)\frac{(\eln k)^2}{\eln \eln k}}\\
&=  e^{\eln k (1 + (1 + \varep^\prime)\frac{\eln k}
{\eln \eln k})
}\\
&< e^{\eln k(1 + \varep)\frac{\eln k}{\eln \eln k}}\\
&= k^{(1 + \varep) \frac{\eln k}{\eln \eln k}}.
\endalign
$$
Here $\varep = \varep^\prime + \varep_3 = 3\varep_3 + \varep_3^{\sssize{2}}$.
  We use (3.6) for the last inequality.

\medskip

Lemmas 17 (iii),  16 (iv), and (3.5) give
$$
\align
|A\{1\}| < k^J k = k^{J+1} &< k^{(1 + \varep_3)\frac{\eln k}{\eln \eln k} + 1} \tag 3.9\\
&< k^{(1 + 2\varep_3)\frac{\eln k}{\eln \eln k}}.
\endalign
$$
The last inequality is again by (3.6).
\medskip
Putting (3.8) and (3.9) together, we have $g(A) < 2k^{(1 + \varep)\frac{\eln k}{\eln \eln k}}$. 
\enddemo

\AuthorRefNames [XXX]
\references
[Bi] \name{Y.\ Bilu}, Structure of sets with small sumset, in
{\it Structure Theory of Set Addition\/}, {\it Ast\'erisque\/} {\bf 258}
(1999), 77--108.

[C1] \name{M.-C.\ Chang}, A polynomial bound in Freiman's theorem,
{\it Duke Math.\ J\/}.\ {\bf 113} (2002), 399--419.

[C2]  \bibline, Factorization in generalized arithmetic
progressions and applications to the Erd\"os-Szemer\'edi sum-product
problems,
{\it GAFA\/}, to appear.

[El] \name{G.\ Elekes}, 
On the number of sums and products,
{\it  Acta Arith\/}.\ {\bf 81}  (1997), 365--367.

[El-R] \name{G.\ Elekes} and \name{I.\  Ruzsa}, Product sets are very large if
sumsets are very small, preprint.

[E] \name{P.\ Erd\H os}, 
  Problems and results on combinatorial number theory.\  {\rm III},
  in {\it Number Theory Day\/}, 43--72
 ({\it Proc.\ Conf. Rockefeller Univ\/}., {\it New York\/}, 1976), 
 {\it Lecture Notes in Math\/}.\ {\bf 626}, Springer-Verlag, New York,
 1977.
 
[E-S] \name{P.\ Erd\H os} and \name{E.\ Szemer\'edi}, 
 On sums and products of integers, {\it Studies in Pure Mathematics},
  Birkh\"auser, Basel,
 1983,
 213--218.

[Fr1] 
\name{G.~A.\ Freiman},
{\it Foundations of a Structural Theory of Set Addition\/},
{\it  Transl.\  of Math.\ Monographs\/} {\bf 37}, A.\ M.\ S.,
Providence, RI, 1973.

[Fr2] \bibline, 
 On the addition of finite sets.\  {\rm I},
{\it   Izv.\ Vysh.~ Ucheb.\ Zaved.\ Matematika\/} {\bf 13} (1959),
202--213.

[Fr3] \bibline,  
Inverse problems of additive number theory.\  {\rm VI}. On the addition of 
finite sets.\  {\rm III},
{\it Izv.~Vysh.~Ucheb.\ Zaved.\ Matematika\/} {\bf 28} (1962), 151--187.
 
[H-T] 
\name{R.\ R.\  Hall} and \name{G.~Tenenbaum},
{\it Divisors\/}, {\it Cambridge Tracts in Math\/}.\ {\bf 90},
Cambridge Univ.\ Press, Cambridge, 1988.

[K-T] 
\name{N.\ Katz} and \name{T.\ Tao},
Some connections between  Falconer's distance set conjecture and
sets of Furstenberg type,
 {\it New York J.\ Math\/}.\ {\bf 7} (2001), 149--157.

[L-R]
\name{M.\ Laczkovich} and \name{I.\ Z.\ Ruzsa}, The number of homothetic subsets, in
{\it The Mathematics of P.\ Erd\H os\/}, {\it II\/} (R.\ L.\ Graham and J.\
Nesetril, eds.), Springer-Verlag, New York, 1977, 294--302.

[Na1] 
\name{M.\ B.~Nathanson},
{\it   Additive Number Theory\/}.\
{\it Inverse Problems and the Geometry of Sumsets\/},
 {\it Grad.\ Texts in Math\/}.\ {\bf 165}, Springer-Verlag, New York,  
1996.
 
[Na2] \bibline, The simplest inverse problems in 
additive number theory, in {\it Number Theory with an Emphasis on the
Markloff Spectrum\/} (Provo, UT, 1991), 191--206,
Marcel Dekker, New York, 1993.
 
[Na3] 
\bibline,   On sums and products of integers, {\it Proc.\
Amer.\ Math.\ Soc\/}.\ {\bf 125} (1997), 9--16.

[N-T]
\name{M.~Nathanson} and \name{G.~Tenenbaum},
 Inverse theorems and the number of sums and products,
in {\it Structure Theory of Set Addition\/},  {\it Ast\'erisque\/} {\bf
258} (1999), 195--204.

[P] 
\name{H. Pl\"unnecke},
 Eine zahlentheoretische Anwendung der Graphtheorie,
{\it J.\ Reine\break Angew.\ Math\/}.\ {\bf 243}
(1970), 171--183.
 
[R] 
\name{W.\ Rudin},
Trigonometric series with gaps,
{\it J.\ Math.\ Mech\/}.\ {\bf 9} (1960), 203--227.

[Ru] 
\name{I.\ Z.\ Ruzsa},
Generalized arithmetical progressions and sumsets,
{\it Acta Math.\ Hungar\/}.\ {\bf 65} (1994), 
379--388.

[Ru2] \bibline,  
Sums of finite sets,
in {\it  Number Theory\/} (New York, 1991--1995), 
 Springer-Verlag, New York, 1996.

[S-T] 
\name{E.~Szemer\'edi} and \name{W.~Trotter},
 Extremal problems in discrete geometry,
{\it Combinatorica\/} {\bf 3} (1983), 381--392.
 
[T] \name{T.\ Tao}, 
From rotating needles to stability of waves: emerging connections
 between combinatorics, analysis, and PDE, {\it Notices Amer.\ Math.\
 Soc\/}.\ {\bf 48} (2001), 294--303.
\endreferences

 \enddocument